\documentclass[12pt]{article}
\usepackage{amsmath,amssymb,latexsym,color,amsfonts,pdfsync,amsthm}
\usepackage{url,graphicx}

\newtheorem{theorem}{Theorem}[section]
\newtheorem*{maintheorem}{Main Theorem}
\newtheorem{lemma}[theorem]{Lemma}

\newtheorem{corollary}[theorem]{Corollary}
\newtheorem{remark}[theorem]{Remark}
\newtheorem{definition}[theorem]{Definition}

\newcommand{\Proof}{ \noindent{\bf Proof:}\quad }

\newcommand{\npmatrix}[1]{\left( \begin{matrix} #1 \end{matrix} \right)}

\def\F{\mathbb{F}}
\def\Fq{\mathbb{F}_q}
\def\PG{\mathrm{PG}}
\def\PGL{\mathrm{PGL}}

\def\GL{\mathrm{GL}}
\def\rk{\mathrm{rk}}

\def\QQ{\mathbb{Q}}
\def\RR{\mathbb{R}}

\def\CC{\mathbb{C}}

\begin{document}

\title{Canonical forms of $2\times 3 \times 3$ tensors over the real field, algebraically closed fields, and finite fields.}
\author{Michel Lavrauw\footnote{The research of the first author was supported by the Fund for Scientific Research - Flanders (FWO) and by a Progetto di Ateneo from Universit\`a di Padova (CPDA113797/11).}
 and John Sheekey\footnote{The second author acknowledges the support of the Fund for Scientific Research - Flanders (FWO).}
}
\date{\today}
\maketitle

\begin{abstract}
We classify the orbits of elements of the tensor product spaces $\F^2\otimes \F^3 \otimes \F^3$ for all finite; real; and algebraically closed fields under the action of two natural groups. The result can also interpreted as the classification of the orbits in the $17$-dimensional projective space of the Segre variety product of a projective line and two projective planes. This extends the classification of the orbits in the $7$-dimensional projective space of the Segre variety product of three projective lines \cite{LaSh2014}. The proof is geometric in nature, relies on properties of the Segre embedding, and uses the terminology of projective spaces. 
\end{abstract}


\section{Introduction}

\subsection{Motivation and main results}

In this article we study orbits of elements of tensor product spaces $V_1\otimes V_2 \otimes V_3$ under the action of two natural groups.
This problem can be studied either in the vector space $V = V_1\otimes V_2 \otimes V_3$ or the projective space $\PG(V_1\otimes V_2 \otimes V_3)$. In this article we will present results using vector space terminology. However we will utilise some terminology and concepts from projective geometry in the proofs.

We consider orbits under the action of the groups $G$ and $H$, where $G$ is the stabiliser in $\GL(V)$ of the set of {\it fundamental tensors} (see Section \ref{sec:prelim} for definitions), and $H$ is the subgroup $\GL(V_1)\times \GL(V_2)\times \GL(V_3)$ of $G$.

These orbits are interesting not only for their own sake, being as they are fundamental objects in both algebra and geometry. They also have many important applications: see for example \cite[Chapter 1]{Landsberg2012}. Furthermore, tensors over finite fields have applications in the study of {\em finite semifields} (that is, nonassociative division algebras), an important and active research topic in finite geometry, see e.g. \cite{LaPo2011}. A deeper understanding of the geometry of three-fold tensor product spaces is crucial to a better structural insight into the set of the isotopism classes of finite semifields, see \cite{Lavrauw2012}. This paper provides a further step in this direction by classifying the orbits in $\F^2\otimes \F^2 \otimes \F^3$ and $\F^2\otimes \F^3 \otimes \F^3$, in the case that $\F$ is a finite field, the real field, or algebraically closed.  The proof uses ideas from \cite{Lavrauw2012}, \cite{LaSh2014}, and \cite{LaSh201?}, and gives geometric characterisations of each orbit.

\begin{maintheorem}
If $\F$ is a finite field, then there are precisely $21$ $H$-orbits and $18$ $G$-orbits of tensors in $\F^2 \otimes \F^3 \otimes \F^3$. For any algebraically closed field $\F$, there are precisely $18$ $H$-orbits and $15$ $G$-orbits of tensors in $\F^2 \otimes \F^3 \otimes \F^3$. There are precisely $20$ $H$-orbits and $17$ $G$-orbits of tensors in $\RR^2 \otimes \RR^3 \otimes \RR^3$.
\end{maintheorem}

This result can be used to fully classify all orbits of subspaces of $\PG(\F^2 \otimes \F^3)$, and hence all orbits of tensors in $\F^2 \otimes \F^3 \otimes \F^r$ for all positive integers $r$. This will be treated by the authors in the forthcoming paper \cite{LaSh23r}.

\subsection{Historical context}

The orbits under the actions of $G$ and $H$ have been studied mostly over algebraically closed fields, see for example \cite{Landsberg2012}. The classification of $G$- (resp. $H$-)orbits of $3\times 3 \times 3$ tensors over the complex numbers has been done in \cite{ThCh1938} (resp. \cite{Nurmiev2000}), and the result for complex $2\times 3 \times 3$ tensors can be obtained from these. However the techniques employed therein (classification of complex ternary cubics in \cite{ThCh1938}, topology in \cite{Nurmiev2000}) are not applicable over finite fields.

The classification of $2\times m \times n$ tensors can be approached using the Weierstrass-Kronecker theory of pencils (that is, $m \times n$ matrices with polynomial entries linear in an indeterminate $\lambda$), see for example \cite{Gantmacher}. From each tensor we can obtain a pencil in a canonical way. However, the classification of pencils does not imply a classification of $G$- or $H$-orbits, as equivalent tensors can give rise to inequivalent pencils. Furthermore, this approach does not give any geometric insight to the nature of these orbits, and does not extend well to larger order tensors.

This approach was used in \cite{Thrall1938} to give canonical forms for $2\times 3 \times 3$ tensors over fields of odd prime order $p>3$. However the forms are presented without proof. For example, there is no proof in the literature that the canonical forms $M_2(x)$ and $M_8(x)$ in \cite[Section 7]{Thrall1938} are indeed the unique orbits with the relevant invariant factors (and in fact these forms are not valid over fields of even characteristic).

Computational and theoretical results over small fields and small dimensions have been obtained, for example \cite{BrHu2012}, \cite{BrSt2012}, \cite{ShGoHa2012}. A geometric classification for the case $\F^2\otimes \F^2 \otimes \F^2$ was given by the authors in \cite{LaSh2014}.

The main result of this paper is to classify orbits of $2\times 3 \times 3$ tensors over all finite fields, and geometrically characterise each orbit. As our approach will be elementary and mostly field-independent, we will also re-obtain the classification over the complex numbers and prime fields, and moreover over any algebraically closed field, and the real numbers.

\section{Definitions and preliminary results}
\label{sec:prelim}

Let $V=\bigotimes_{i=1}^r V_i$, where $V_1, \ldots, V_r$ are finite dimensional vector spaces over some field $\F$, with $\dim V_i=n_i<\infty$. The set of {\it fundamental tensors} is the set $\{v_1\otimes v_2 \otimes \ldots \otimes v_r :  v_i\in V_i\backslash\{0\}\}$. Projectively, this corresponds to the {\it Segre variety} $S_{n_1,n_2,\ldots,n_r}(\F)$, that is the image of a {\it Segre embedding} $\sigma_{n_1,\ldots,n_r}$ defined as
\begin{align*}
\sigma_{n_1,\ldots,n_r}&:\PG(V_1)\times \PG(V_2)\times \ldots \PG(V_r) ~\rightarrow ~\PG(V)\\
&:(\langle v_1\rangle ,\langle v_2\rangle ,\ldots, \langle v_r\rangle )\mapsto 
\langle v_1\otimes v_2 \otimes \ldots \otimes v_r\rangle.
\end{align*}


The {\it rank} of a tensor $A$, denoted $\rk(A)$, is defined to be the minimum number $k$ such that there exist fundamental tensors $\alpha_1,\ldots,\alpha_k$ with $A \in \langle \alpha_1,\ldots,\alpha_k\rangle$. The fundamental tensors hence correspond precisely to the set of {\it rank one tensors}. When $r=2$, that is $V = V_1 \otimes V_2$, the tensor rank corresponds to the usual matrix rank. We define the rank of a projective point in $\PG(V)$ to be the rank of any corresponding tensor.

The setwise stabilizer in $\GL(V)$ of the set of rank one tensors in $V = V_1\otimes V_2 \otimes V_3$ will be denoted by $G$.  Clearly the rank of a tensor is $G$-invariant. The group $\GL(V_1)\times\GL(V_2)\times\GL(V_3)$ defines a subgroup of $G$, which we will denote by $H$, via the action defined by
\[
(g_1,g_2,g_3): v_1\otimes v_2 \otimes v_3 \mapsto  v_1^{g_1}\otimes v_2^{g_2} \otimes v_3^{g_3}.
\]
Note that $H$ is not necessarily equal to $G$, for example if $V_1 \ne V_2=V_3$, then $G\cong \GL(V_1) \times (\GL(V_2)\wr \mathrm{Sym}(2))$. The group $\mathrm{Sym}(2)$ arises from the map $T$ defined by $T:v_1\otimes v_2 \otimes v_3 \mapsto  v_1\otimes v_3 \otimes v_2$.

For $A\in V_1\otimes V_2 \otimes V_3$ we define the {\it first contraction space of $A$} as the following subspace of
$V_2\otimes V_3$:
$$A_1:=\langle w_1^\vee(A)~:~w_1^\vee \in V_1^\vee\rangle,$$
where $V_1^\vee$ denotes the dual space of $V_1$, and the {\it contraction} 
$w_1^\vee(A)$ is defined by its action on the fundamental tensors 
\begin{eqnarray}
w_1^\vee(v_1\otimes v_2 \otimes v_3)=w_1^\vee(v_1)v_2\otimes v_3.
\end{eqnarray}
Similarly we define the {\it second} and {\it third contraction space}, and denote these 
by $A_2$ and $A_3$ respectively. We will consider the projective subspaces $\PG(A_i)$ of $\PG(V_j \otimes V_k)$, where $j<k$ and $\{i,j,k\}=\{1,2,3\}$.

\bigskip



The setwise stabilizer of the set of rank one tensors in the contracted space $V_j \otimes V_k$ will be denoted by $G_i$, where $j<k$ and $\{i,j,k\}=\{1,2,3\}$, and the subgroups $\GL(V_j)\times\GL(V_k)$ of $G_i$ by $H_i$. For example, with $V_1=\F^2$, $V_2=V_3=\F^n$, $n>2$, the group $G_1\cong \GL(n,q)\wr Sym(2) > H_1$, while $G_2=H_2$ and $G_3=H_3$.

The following elementary lemmas will prove useful in the classification to follow.

\begin{lemma}\label{lemma:contraction_equivalence}
If $A, B\in V_1\otimes V_2 \otimes V_3$ then the following are equivalent:
\begin{enumerate}
\item[(i)] $A\sim_H B$;
\item[(ii)] $\forall i$ $A_i \sim_{H_i} B_i$;
\item[(iii)] $\exists i$  $A_i \sim_{H_i} B_i$.
\end{enumerate}
\end{lemma}
\Proof
Suppose $A_1 \sim_{H_1} B_1$, i.e. suppose there exists $(h_2,h_3)\in H_1$ such that
$A_1^{(h_2,h_3)}=B_1$. Let $A=\sum a_{ijk} u_i\otimes v_j \otimes w_k$, where the $u_i$'s, $v_i$'s, $w_i$'s form a basis of $V_1$, $V_2$, $V_3$, respectively. Then
\begin{eqnarray}
A_1=\langle \sum a_{ijk}v_j\otimes w_k~:~i=1,...,\dim V_1\rangle, 
\end{eqnarray}
and 
\begin{eqnarray}
B_1=A_1^{(h_2,h_3)}=\langle \sum a_{ijk}v_j^{h_2}\otimes w_k^{h_3}~:~i=1,...,\dim V_1\rangle.
\end{eqnarray}
This implies that $B=\sum b_iu'_i\otimes (a_{ijk} v_j^{h_2}\otimes w_k^{h_3})$ where $b_i\in \F$, and the $u'_i$'s form a basis. Hence
$B=\sum u_i^{h_1}\otimes (a_{ijk} v_j^{h_2}\otimes w_k^{h_3})=A^{(h_1,h_2,h_3)}$, where $\GL(V_1)\ni h_1$ maps $u_i$ to $b_i u'_i$ for $b_i\neq 0$ and $u_i$ to $u'_i$ otherwise.
This proves $(iii) \Rightarrow (i)$. The proof of  $(i) \Rightarrow (ii)$ is similar, and  $(ii) \Rightarrow (iii)$ is trivial.
\qed

\begin{corollary}\label{lemma:contraction_equivalence_G}
Suppose $V_1 \ne V_2=V_3$. Then two tensors $A$ and $B$ are $G$-equivalent if and only if $A_1$ is $H_1$-equivalent to one of $\{B_1,B^T_1\}$, if and only if $A_1$ is $G_1$-equivalent to $B_1$.
\end{corollary}

Suppose $V_2=V_3 = \F^n$, and let $\{e_1,\ldots,e_n\}$ be the standard basis for $\F^n$. Then $V_2 \otimes V_3 \simeq M_n(\F)$, and we can choose the isomorphism such that the tensor $x := \sum_{i=1}^n e_i \otimes e_i$ corresponds to the identity matrix. Then the stabiliser of $x$ in $H_1$ is given by $N := \{(g,g^{-1}):g \in \GL(n,\F)\}$.

\begin{lemma} \label{lem:perp}
If $S_U$ denotes the stabiliser in $\mathrm{GL}(V)$ of a $k$-dimensional subspace $U$ of a vector space $V$. Then $\{g^{-1}:g \in S_U\}$ is the stabiliser in ${\mathrm{GL}}(V)$ of a subspace of codimension $k$.
\end{lemma}
\Proof
Let $V_1$ denote the $n$-dimensional vector space of row vectors and $V_2$ the $n$-dimensional vector space of column vectors over the field $\F$. Then $\mathrm{GL}(V_1)=\mathrm{GL}(V_2)$ with the action from the right on $V_1$ and from the left on $V_2$. Let $U$ be a subspace of $V_1$ and define $U^\perp:=\{y\in V_2~:~ \forall u \in U~:~uy=0\}$. Suppose $A\in \mathrm{GL}(V_1)$ is in the stabiliser $S_U$ in ${\mathrm{GL}}(V_1)$ of $U$. Then $\forall u \in U$, $uA^{-1}\in U$ and hence $uA^{-1}y=0$, $\forall y\in U^\perp$. But this means that $A^{-1}U^\perp \subseteq U^\perp$. It follows that $\{A^{-1}~:~ A \in S_U\}$ is the stabiliser of $U^\perp$ in ${\mathrm{GL}}(V_2)$.
\qed


\begin{lemma}\label{lem:Q(x)}
A rank two point $x$ in $\langle S_{n_1,n_2}(\F)\rangle$ is contained in a unique solid which intersects
$S_{n_1,n_2}(\F)$ in a subvariety $Q(x)$ equivalent to a Segre variety $S_{2,2}(\F)$.
\end{lemma}
\Proof
Consider a point $x$ of rank two, i.e. a point $x\in \langle y,z\rangle$ where $y,z \in S_{n_1,n_2}(\F)$, $y\neq x\neq z$. If $y=\sigma_{n_1,n_2}(y_1\times y_2)$ and $z=\sigma_{n_1,n_2}(z_1\times z_2)$, then
$Q_{y,z}:=\sigma_{n_1,n_2}(\langle y_1,z_1\rangle \times \langle y_2,z_2\rangle)$ is a Segre variety equivalent to $S_{2,2}(\F)$ contained in $S_{n_1,n_2}(\F)$, and $x\in \langle Q_{y,z}\rangle$. Any subvariety of $S_{n_1,n_2}(\F)$ equivalent to $S_{2,2}(\F)$ is of the form $Q_{y,z}$ for some
$y,z \in S_{n_1,n_2}(\F)$, and two solids spanned by such subvarieties 
$Q_{y,z}=\sigma_{n_1,n_2}(\langle y_1,z_1\rangle \times \langle y_2,z_2\rangle)$ and
$Q_{y',z'}=\sigma_{n_1,n_2}(\langle y'_1,z'_1\rangle \times \langle y'_2,z'_2\rangle)$ intersect in the subspace
spanned by the subvariety $\sigma_{n_1,n_2}(w_1\times w_2)$, where
$w_1 =\langle y_1,z_1\rangle \cap \langle y'_1,z'_1\rangle$ and $w_2= \langle y_2,z_2\rangle \cap \langle  y'_2,z'_2\rangle$. It follows that the intersection of two different solids $\langle Q_{y,z}\rangle$ and $\langle Q_{y',z'}\rangle$ is either a point on $S_{n_1,n_2}(\F)$ or a line on $S_{n_1,n_2}(\F)$, and hence cannot contain a point of rank two.
This implies that $Q(x):= Q_{y,z}$ is uniquely defined by $x$.
\qed

\begin{definition}
The {\it ($i$-th) rank distribution} $r_i(A)$ of a tensor $A$ is defined to be the tuple whose $j$-th entry is the number of rank $j$ points in the $i$-th contraction space $\PG(A_i)$.
\end{definition}
Note that if $A$ and $B$ are $H$-equivalent, then $r_i(A)=r_i(B)$ for each $i$. The converse does not necessarily hold. Also note that the corresponding statement for $G$-equivalence is false.

\section{Orbits of tensors in $\F^2\otimes \F^3 \otimes \F^3$}

In this section we study the orbits of tensors in $\F^2 \otimes \F^3 \otimes \F^3$ in terms of their first contraction spaces, by considering the corresponding projective subspaces of $\PG(\F^3 \otimes \F^3) \simeq \PG(8,\F)$. By Lemma \ref{lemma:contraction_equivalence} and its corollary, the $H_1$-orbits (resp. $G_1$-orbits) of these subspaces are in one-to-one correspondence with the  $H$-orbits (resp. $G$-orbits) of nonzero tensors. The proof itself contains a geometric characterisation of each of the orbits. A nonzero tensor in $\F^3 \otimes \F^3$ (and similarly a point in $\PG(\F^3 \otimes \F^3)$) is called {\it singular} if it has rank at most two.

\subsection{The trivial cases $\mathbf{o_0,o_1,o_2,o_3}$}\label{sec:trivial}
The zero vector forms an $H$-orbit, which we denote by $\mathbf{o_0}$. 
The next trivial cases occur when $\PG(A_1)$ is a point, and this gives rise
 to three further orbits: $\mathbf{o_1,o_2,o_3}$. Namely, if $A\in \F^2\otimes \F^3 \otimes \F^3$, $A\neq 0$, with first contraction space $A_1$ of dimension one, then $A$ is of the form
$x\otimes \tau$, with $x\in \F^2$ and $\tau \in \F^3 \otimes \F^3$. Two such nonzero tensors
$A=x\otimes \tau$ and $A'=x'\otimes \tau'$ are $G$-equivalent if and only if $\tau$ and $\tau'$ are
$G_1$-equivalent. Therefore, it is clear that
there are three orbits on the nonzero vectors in 
$\{A\in \F^2\otimes \F^3 \otimes \F^3~:~\mbox{dim}(A_1)=1\}$. Call these orbits $\mathbf{o_1}$, $\mathbf{o_2}$ and $\mathbf{o_3}$, where
the subscript $i$ corresponds to the rank of the 2-fold tensor, i.e. $A\in o_i$ if and only if $\rk(A_1)=i$.

\subsection{Lines containing at least one singular point}\label{sec:singular}

Next we consider tensors $A\in\F^2\otimes \F^3 \otimes \F^3$ whose first contraction space $A_1$ is two-dimensional and contains a singular tensor. This means that $A_1$ determines a line $\PG(A_1)$ containing a singular point.

\begin{theorem}
\label{thm:main}
Let $\Omega$ denote the set of tensors  in $\F^2 \otimes \F^3 \otimes \F^3$ whose first contraction space is two-dimensional and contains at least one singular tensor. If $\F$ is finite, or $\F$ is real closed, there are $16$ $H$-orbits and $13$ $G$-orbits on $\Omega$. If $\F$ is algebraically closed, there are $14$ $H$-orbits and $11$ $G$-orbits on $\Omega$.
\end{theorem}

\Proof 
We will consider the orbits of points and lines in $\PG(\F^3 \otimes \F^3)$ under the action induced by $H_1$. We divide the proof into subcases by considering the rank distribution. For a tensor $A \in \F^2\otimes \F^3\otimes \F^3$ we denote the rank distribution $r_1(A)$ of $A$ by $[a,b,c]$. 

\medskip

Let $p_1, p_2, p_3, p_4$ be the standard frame of $\PG(\F^3)$ corresponding to the basis $e_1,e_2, e_3$ of $\F^3$, and recall the Segre embedding
\begin{eqnarray}
\sigma_{3,3}~:~\PG(\F^3)\times \PG(\F^3)\rightarrow S_{3,3}(\F)\subset \PG(\F^3\otimes \F^3).
\end{eqnarray}

\medskip


{\bf Case (1):} {\it $\PG(A_1)$ is a line, with $a \geq 2$; $\mathbf{o_4,o_4^T,o_5}$.}

\medskip

We are supposing $a\geq 2$, i.e. $\PG(A_1)$ has at least two points of rank 1. The lines contained in $S_{3,3}(\F)$ are either of the form $\sigma_{3,3}(p_1\times \langle p_1,p_2\rangle)$ or
$\sigma_{3,3}(\langle p_1,p_2\rangle \times p_1)$
and they form two orbits under the action induced by $H_1$. Clearly these are $G_1$-equivalent under the action induced by the element $T:v_2 \otimes v_3 \mapsto v_3 \otimes v_2$, and so we denote them by $\mathbf{o_4}$ and $\mathbf{o_4^T}$ respectively. In the case $\F=\F_q$, these lines have rank distribution $[q+1,0,0]$.

\medskip

The lines not contained in $S_{3,3}(\F)$ have $a\leq 2$, since $S_{3,3}(\F)$ is the intersection of quadrics. 
Suppose $a=2$. Then $\PG(A_1)$ is a two-secant to $S_{3,3}(\F)$ and hence determines a unique subvariety $S_{2,2}(\F)$.
Each such subvariety is equivalent to $\sigma_{3,3}(\langle p_1,p_2\rangle \times \langle p_1,p_2\rangle)$, and its stabiliser inside the group induced by $H_1$ acts 3-transitively on the points of $\langle p_1,p_2\rangle$ in each factor.
Therefore, the set of two-secants forms one orbit under the group induced by $H_1$. All other points of $\PG(A_1)$ have rank two and hence if $\F=\F_q$, the line $\PG(A_1)$ has rank distribution $[2,q-1,0]$. We denote this orbit by $\mathbf{o_5}$.

\bigskip

{\bf Case (2):} {\it $\PG(A_1)$ is a line, with $a = 1$; $\mathbf{o_6,o_7,o_7^T,o_8,o_9}$.}

\bigskip

Next suppose that $a=1$ and $b\geq 1$, i.e. $\PG(A_1)$ is tangent to $S_{3,3}(\F)$ and contains at least one point of rank 2.
Assume $\PG(A_1)=\langle x_1, x_2\rangle$, with $x_1$ of rank one and $x_2$ or rank two.
We consider three cases depending on the position of $x_1$ with respect to $Q(x_2)$. 

\bigskip

$\mathbf{o_6}$. If $x_1$ is contained in $Q(x_2)$ then $\PG(A_1)$ is a tangent line to $Q(x_2)$ in $\langle Q(x_2)\rangle$ 
and it follows from \cite{LaSh2014} that that these lines form one orbit under the group induced by $H_1$. 
This orbit has rank distribution $[1,q,0]$ in case $\F=\F_q$.

\bigskip

$\mathbf{o_7}$. If $x_1$ is not contained in $Q(x_2)$, but is contained in a subvariety $S_{2,3}(\F)$ of $S_{3,3}(\F)$ containing $Q(x_2)$, then 
w.l.o.g. we may assume that $\PG(A_1)=\langle x_1,x_2\rangle$ is contained in $\PG(\langle e_1,e_2\rangle \otimes\F^3)$, where
$x_2=\langle e_1\otimes e_1 + e_2\otimes e_2\rangle$ and 
$x_1=\langle e_1\otimes e_3\rangle$, since the subgroup
of the group induced by $H_1$ stabilising $\sigma_{3,3}(\langle p_1,p_2\rangle \times \PG(\F^3))$ acts transitively on frames of both factors.
Note that the dimension
of the third contraction space $A_3$ is three, since otherwise we would be in the previous case $o_6$.
In the finite case $\F=\F_q$, the line $\PG(A_1)$ has rank distribution $[1,q,0]$. 

\bigskip

$\mathbf{o_7^T}$. If $x_1$ is not contained in $Q(x_2)$, but is contained in a subvariety $S_{3,2}(\F)$ of $S_{3,3}(\F)$ containing $Q(x_2)$, repeating the argument for $o_7$ we see that we get a single orbit, which is $G_1$-equivalent to $o_7$, and also has rank distribution $[1,q,0]$. Note that $o_7$ and $o_7^T$ cannot be $H_1$-equivalent, for then $\PG(A_1)$ would be contained in some $\langle S_{2,2}(\F)\rangle$.


\bigskip

$\mathbf{o_8}$. If $x_1$ is not contained in a subvariety $S_{2,3}(\F)$ or $S_{3,2}(\F)$ of $S_{3,3}(\F)$ containing $Q(x_2)$, 
then, again using the transitivity of the group induced by $H_1$ on frames of both factors, w.l.o.g. we have
$x_1=\langle e_1\otimes e_1\rangle$ and $x_2=\langle e_2\otimes e_2 + e_3\otimes e_3\rangle$.
If $\F=\F_q$ then $\PG(A_1)$ has rank distribution $[1,1,q-1]$.

\bigskip

$\mathbf{o_9}$. This leaves only one more possibility for tangent lines, namely when $b=0$ and $a=1$. An example of such a line
is $\langle x_1,x_2\rangle$ with $x_1=\langle e_1\otimes e_3\rangle$ and $x_2=\langle e_1\otimes e_1 + e_2\otimes e_2 + e_3\otimes e_3\rangle$. We show that also these lines form one orbit. For suppose $\PG(A_1)$ is a line containing
a point $x_1$ of rank one and a point $x_2$ of rank 3. Then w.l.o.g. we may assume that
$x_2=\langle e_1\otimes e_1 + e_2\otimes e_2 + e_3\otimes e_3\rangle$ and 
$x_1=\langle x'\otimes x''\rangle $ with $x'=ue_1+ve_2+we_3$ for some $u,v,w \in \F$ with $u\neq 0$. Now $x_2$ can also be written as
$$
x_2=\langle (ue_1+ve_2+we_3)\otimes \frac{e_1}{u} + e_2\otimes (e_2-\frac{v}{u}e_1)+e_3\otimes (e_3-\frac{w}{u}e_1)\rangle.
$$
After an appropriate base change we see that up to $H_1$-equivalence we may assume
$x_1=e_1\otimes (u'e_1+v'e_2+w'e_3)$, for some $u',v',w' \in \F$ with $w'\neq 0$, and $x_2=e_1\otimes e_1+e_2\otimes e_2 + e_3\otimes e_3$. Now
note that necessarily $u'=0$, since otherwise the line $\langle x_1,x_2\rangle$ would contain the rank two 
point $\langle x_2-u'x_1\rangle$ in contradiction with our hypothesis $b=0$. So w.l.o.g. we have
$x_1=e_1\otimes (v'e_2+w'e_3)$.
Finally, since $x_2$ can be rewritten as
$$
e_1\otimes e_1+ (e_2-\frac{v'}{w'}e_3)\otimes e_2 + \frac{e_3}{w'}\otimes(v'e_2+w'e_3),
$$
we see that the line $\langle x_1,x_2\rangle$ is equivalent to 
$\langle e_1\otimes e_3, e_1\otimes e_1 + e_2\otimes e_2 + e_3\otimes e_3\rangle$, as claimed.
We denote this orbit by $\mathbf{o_9}$. If $\F=\F_q$ these lines have rank distribution $[1,0,q]$.

\bigskip

The orbits $o_1, \ldots, o_9$ are the only orbits on lines with rank distribution $[a,b,c]$ with $a\geq 1$. So from now on we assume $a=0$.

\bigskip

{\bf Case (3):} {\it $\PG(A_1)$ is a line, with $a=0$, $b\geq 2$; $\mathbf{o_{10},o_{11},o^T_{11},o_{12},o_{13},o_{14}}$.}

\bigskip

$\mathbf{o_{10}}$. Our next orbit is inherited from the $2\times 2\times 2$ case; it consists of lines which are contained in a three-dimensional subspace spanned by a subvariety $S_{2,2}(\F)$ of $S_{3,3}(\F)$ but disjoint from this subvariety. These corresponds to the orbit of nonsingular points in the Segre variety product of three projective lines, see \cite{LaSh2014}. When $\F=\Fq$ or $\F=\RR$, by \cite{LaSh2014} these lie in one orbit which we denote by $\mathbf{o_{10}}$. The rank distribution
for $\F=\F_q$ is $[0,q+1,0]$, i.e. {\it constant rank two}. Note that this orbit is empty if $\F$ is algebraically closed.

$\mathbf{o_{11}}$. Another orbit containing constant rank two lines, consists of the lines disjoint from $S_{3,3}(\F)$ which are contained in a five-dimensional space spanned by some subvariety $S_{2,3}(\F)$, but not any $S_{2,2}(\F)$. We claim that all these lines are equivalent to the line $\langle e_1\otimes e_1+e_2\otimes e_2,e_1\otimes e_2+e_2\otimes e_3\rangle$. 
Suppose $\langle x_1,x_2\rangle$ is a constant rank two line contained in a subvariety $S_{2,3}(\F)$. Then w.l.o.g. we 
may assume that $\langle x_1,x_2\rangle$ is contained in $\PG(\langle e_1,e_2\rangle \otimes \F^3)$. It follows that
$Q(x_1)$ and $Q(x_2)$ meet in a line $\ell$, and again w.l.o.g. we may assume 
$\ell=\sigma_{3,3}(\langle e_1,e_2\rangle \otimes \langle e_2\rangle)$. Each plane $\pi_i=\langle x_i,\ell\rangle$, $i=1,2$, intersects the quadric $Q(x_i)$ in two lines: the line $\ell$ and another line, say $\ell_i$. If $\ell_1$ and $\ell_2$ have a point in common, then the plane $\langle \ell_1,\ell_2\rangle$ is contained in the Segre variety $S_{3,3}(\F)$. In this case the three-dimensional space $\langle \ell, \ell_1,\ell_2\rangle$ contains the plane $\langle \ell_1,\ell_2\rangle$ as well as the line $\langle x_1,x_2\rangle$. But then $\langle x_1,x_2\rangle$ meets $\langle \ell_1,\ell_2\rangle$ contradicting our hypothesis that $\langle x_1,x_2\rangle$ is a constant rank two line. Therefore $\ell_1\cap \ell_2=\emptyset$. W.l.o.g. we may assume that $y_1=\ell\cap \ell_1=\langle e_1\otimes e_2\rangle$ and $y_2=\ell\cap \ell_2=\langle e_2\otimes e_2\rangle$. Now considering the intersection point $\langle e_1 \otimes \alpha \rangle$ of the line $\ell_1$ with the line $\langle x_1,y_2\rangle$ and the intersection point $\langle e_2\otimes \beta \rangle$ of the line $\ell_2$ with the line  $\langle x_2,y_1\rangle$ we have that $x_1=\langle ue_1\otimes \alpha+e_2\otimes e_2\rangle$ and 
$x_2=\langle e_1\otimes e_2+ve_2\otimes \beta \rangle$, $u\neq 0\neq v$. Since the $\alpha, e_2, \beta$ are 
linearly independent, we may assume $u\alpha=e_1$ and $v\beta=e_3$. This proves the claim.
We call this orbit $\mathbf{o_{11}}$.

\bigskip

$\mathbf{o_{11}^T}$. Replacing $S_{2,3}(\F)$ with $S_{3,2}(\F)$ in the above, and repeating the argument, we get another orbit $\mathbf{o_{11}^T}$, which is $G_1$-equivalent to $o_{11}$. These can not be $H_1$-equivalent, for if they were then $\PG(A_1)$ would be contained in some $\langle S_{2,2}(\F)\rangle$, contradicting our assumption.

\bigskip

The next step in the proof is to consider the remaining orbits of lines which have at least two points of rank two, and no points of rank one.
The assumption that such a line $\langle x_1,x_2\rangle$, where $x_1$ and $x_2$ are points of rank two, is not in one of the previous orbits implies that it is not contained in any subspace
spanned by a subvariety $S_{2,3}(\F)$ or $S_{3,2}(\F)$ of $S_{3,3}(\F)$, and hence $Q(x_1)$ and $Q(x_2)$ must intersect in a point, say $z$.
The orbits of such lines depend on the position of the points $x_1,x_2,z$. We have the following three possibilities:
(i) none of $x_1$, $x_2$ lies on a secant through $z$;  (ii) only one of $x_1$, $x_2$ lies on a secant throught $z$; (iii) both $x_1$ and $x_2$ lie on a secant through $z$.

\bigskip
(i) $\mathbf{o_{12}}$. In this case we assume w.l.o.g. that $z=\langle e_1\otimes e_2\rangle $. Denote the two planes of the Segre variety $S_{3,3}(\F)$ that contain $z$ by $\pi(z)=\PG(\F^3\otimes e_2)$ and $\pi'(z)=\PG(e_1\otimes \F^3)$.
Then the hypotheses imply that for each $i\in \{1,2\}$, $x_i$ lies on the tangent plane of $Q(x_i)$ at $z$, and hence
$x_i$ is on a secant line $\langle y_i,y'_i\rangle$ with $\langle z,y_i\rangle =\pi(z)\cap Q(x_i)$ and 
$\langle z,y'_i\rangle =\pi'(z)\cap Q(x_i)$. Summarizing, w.l.o.g. we may assume that
$x_i=\langle \alpha_i\otimes e_2+e_1\otimes\alpha'_i\rangle$, where
both $e_1,\alpha_1,\alpha_2$ and $e_2,\alpha'_1,\alpha'_2$ are linearly independent. It follows that
the line $\langle x_1,x_2\rangle$ is equivalent to the line
$\langle e_1\otimes e_1+e_2\otimes e_2,e_1 \otimes e_3 + e_3\otimes e_2\rangle$, under the action induced by $(g_1,g_2)$, where $g_1:e_1\mapsto e_1,\alpha_1 \mapsto e_2,\alpha_2\mapsto e_3$, $g_2:e_2\mapsto e_2,\alpha'_1 \mapsto e_1,\alpha'_2\mapsto e_3$. We call this orbit $\mathbf{o_{12}}$, and observe that a line in this orbit is a constant rank two line not contained in the span of any subvariety $S_{2,3}(\F)$ or $S_{3,2}(\F)$.

\bigskip

(ii) $\mathbf{o_{13}}$. In this case we assume w.l.o.g. that $z=\langle e_1\otimes e_2\rangle$, and that $x_2$ lies on a secant through $z$.
Then by hypothesis $x_1$ does not lie on a secant through $z$, hence, using the same notation as in 
case (i) $x_1$ lies on the tangent plane of $Q(x_1)$ at $z$, and therefore on a secant line $\langle y_1,y'_1\rangle$ with $\langle z,y_1\rangle =\pi(z)\cap Q(x_1)$ and 
$\langle z,y'_1\rangle =\pi'(z)\cap Q(x_1)$. Continuing the argument as before, we may conclude that w.l.o.g.
$x_1=\langle e_1\otimes e_1 +e_2\otimes e_2\rangle $ and $x_2=\langle e_1\otimes e_2+ e_3\otimes e_3\rangle$.
We denote this orbit by $\mathbf{o_{13}}$; its lines have rank distribution $[0,2,q-1]$  over $\F_q$.

\bigskip

(iii) $\mathbf{o_{14}}$. In this case w.l.o.g. we assume $z=e_2\otimes e_2$. Then there are $\alpha,\beta,\gamma,\delta$ such that $x_1=\alpha\otimes\beta +e_2\otimes e_2$ and 
$x_2=e_2\otimes e_2+\gamma\otimes\delta$. Since $\langle x_1,x_2\rangle$ is not contained in any subspace
spanned by a subvariety $S_{2,3}(\F)$ or $S_{3,2}(\F)$ of $S_{3,3}(\F)$, we have that both $\alpha,e_2,\gamma$ and $\beta,e_2,\delta$ are linearly independent. It follows that each such line is equivalent to $\langle x_1,x_2\rangle$ with
$x_1=e_1\otimes e_1 +e_2\otimes e_2$ and 
$x_2=e_2\otimes e_2+ e_3 \otimes e_3$. 
We call this orbit $\mathbf{o_{14}}$; its lines have rank distribution $[0,3,q-2]$ over $\F_{q}$.

\bigskip

We have now dealt with all the orbits on lines which have at least one point of rank one and with all the orbits
on lines with at least two points of rank two. So from now on we only
need to consider the lines $\PG(A_1)$ with rank distribution $[a,b,c]$ where $a=0$ and $b\leq 1$.

\bigskip

{\bf Case (4):} {\it $\PG(A_1)$ is a line, with $a=0$, $b = 1$; $\mathbf{o_{15},o_{16}}$.}

If $b=1$, then let $x_1$ be a point of rank three and $x_2$ a point of rank two on $\PG(A_1)$.
We distinguish the following two cases: a) there exist two points $y_i$, $i=1,2$ of rank $i$ such that
$x_1$ is on the line $\langle y_1,y_2\rangle$ and $Q(y_2)=Q(x_2)$; and b) such points do not exist.

\bigskip

a) $\mathbf{o_{15}}$. In this case w.l.o.g. we may assume that $y_1=\langle e_3\otimes e_3\rangle$ and 
$Q(y_2)= Q(x_2) =\sigma_{3,3}(\langle e_1, e_2\rangle \times \langle e_1,e_2\rangle)$.
Now projecting the line $\langle x_1,x_2\rangle$ from $y_1$ onto $\langle Q(x_2)\rangle$
gives the constant rank two line $\langle x_2,y_2\rangle$ in $\langle Q(x_2)\rangle$ (cf. orbit $o_{10}$). 
Using the transitivity of the stabilizer of $Q(x_2)$ inside the group induced by $H_1$ on these lines (cf. orbit $o_{10}$), we
obtain one orbit of such lines. We denote this orbit by $\mathbf{o_{15}}$. A representative for the orbit is
$\langle e_1\otimes e_1+ e_2\otimes e_2+u e_1\otimes e_2+e_3\otimes e_3, 
e_1\otimes e_2+v e_2\otimes e_1\rangle$, with $u,v \in \F$ such that $v\lambda^2+uv\lambda - 1 \neq 0$ for all
$\lambda \in \F$.

If $\F=\Fq$, the orbit $o_{15}$ has rank distribution $[0,1,q]$. If $\F$ is algebraically closed, this orbit is empty.

\bigskip

b) $\mathbf{o_{16}}$. In this case let $Q(x_2)=\sigma_{3,3}(\ell \times m)$, and suppose 
$x_1=\langle e_1\otimes e_1 + e_2\otimes e_2 + e_3 \otimes e_3 \rangle$.
Then \\
\\
$(**)$ $x_1$ can be written as
$x_1=\langle a_1\otimes b_1 + a_2\otimes b_2 + a_3 \otimes b_3 \rangle$,
if and only if there exists a nonzero constant $k\in \F$ such that 
$a_i \cdot b_j=0$ for $i\neq j$ and $a_i\cdot b_i=k$, where $\cdot$ denotes the usual dot product of vectors.\\
\\
Recall that an element of the stabiliser $N$ of $x_1$ has the form $(g,g^{-1})$, and corresponds to conjugation in ${\mathrm{GL}}(3,\F)$.
By Lemma \ref{lem:perp}, the set of elements $(g,g^{-1})$ in $N$ such that $g$ stabilises a point (line) $x$ in the first factor 
corresponds to the set of elements $(g,g^{-1})$ such that $g^{-1}$ stabilises a line (point)
in the second factor, and vice versa.

Now assume $\ell=\langle a_1,a_2\rangle$, and consider the duality from $\PG(\F^3)$ (the first factor in the
pre-image of $\sigma_{3,3}$) to $\PG(\F^3)$ (the second factor in the pre-image of $\sigma_{3,3}$) 
induced by the standard inner
product. The dual space of a subspace will be denoted by $\perp$. If $m$ does not pass through 
$\langle b_3\rangle:=\ell^\perp$, then put
$\langle b_1\rangle=\langle a_2\rangle^\perp \cap m$, 
$\langle b_2\rangle=\langle a_1\rangle^\perp \cap m$, 
$\langle a_3\rangle = m^\perp$.
After appropriate scaling we have the necessary conditions $(**)$ to write 
$x_1=\langle a_1\otimes b_1 + a_2\otimes b_2 + a_3 \otimes b_3 \rangle$
with 
$Q(x_2)=\sigma_{3,3}(\langle a_1,a_2\rangle \times \langle b_1,b_2\rangle)$.
This is in contradiction with the hypothesis (put $y_2=\langle a_1\otimes b_1+a_2\otimes b_2\rangle$ and $y_1=a_3\otimes b_3 \rangle$).
It follows that $m$ is a line through $\langle b_3\rangle=\ell^\perp$. Let
$m=\langle b_2,b_3\rangle$ and suppose
$x_2=\langle  a_1\otimes (a b_2 + b b_3) + a_2\otimes (c b_2 + d  b_3)\rangle$ for some $a,b,c,d \in \F$. Then
$c=0$ since otherwise $\langle x_1, x_2 \rangle$ would have a rank two point
$\langle a_1\otimes (-cb1+ab_2+bb_3) + (da_2-ca_3)\otimes b_3\rangle$ which is different from $x_2$, a contradiction.
Arguing as before it follows that there is one orbit, which we call $o_{16}$, 
satisfying the hypotheses and it is
represented by the line $\langle x_1,x_2 \rangle$ with
$x_1=\langle e_1\otimes e_1 + e_2\otimes e_2 + e_3 \otimes e_3 \rangle$ and
$x_2=\langle e_1\otimes e_2 + e_2 \otimes e_3\rangle$. If $\F=\Fq$, this orbit has rank distribution $[0,1,q]$, and is denoted by $\mathbf{o_{16}}$.

This completes the classification of lines containing at least one point of rank less than $3$.  \qed

\begin{remark}
The existence of lines in orbits $o_{10}$ and $o_{15}$ depends on the field $\F$. The orbit $o_{10}$ corresponds to the orbit of a nonsingular point in the Segre variety product of three projective lines, see \cite{LaSh2014}. These in turn correspond to algebraic field extensions of degree $2$. If $\F$ is algebraically closed, then there are no such extensions and hence the orbit $o_{10}$ does not occur. The same holds true for the orbit $o_{15}$ (see proof).
For any other field $\F$, if we let $\eta(\F)$ denote the number of isotopism classes (see \cite[Theorem 4.3]{Lavrauw2012}) of algebraic field extensions of degree $2$ of $\F$ (which may be infinite), then the number of $H$-orbits on $\Omega$ is
$14+2 \eta(\F)$.
\end{remark}

\subsection{Lines without singular points}

The remaining case, that is the case of lines of constant rank $3$, is field-dependent, and it will be convenient to utilise the language of matrices, determinants and characteristic polynomials. 

Suppose the line $\langle x_1,x_2\rangle$, $x_i=\langle v_i\rangle$, has rank distribution $[0,0,q+1]$, that is, a constant rank $3$ line. We may identify this with a $2$-dimensional subspace $\langle v_1,v_2\rangle$ of $M_3(\F)$. Then we have that $f(t) := \det( v_1 - t v_2) \ne 0$ for all $t\in \F$. But $f$ is a polynomial of degree $3$ over $\F$, and hence if $\F$ is algebraically closed (or real closed), must have a root in $\F$, a contradiction. So, taking into account the results from Sections \ref{sec:trivial} and \ref{sec:singular}, we have shown the following.
\begin{theorem}
For any algebraically closed (or real closed) field  $\F$, there are precisely $18$ $H$-orbits and $15$ $G$-orbits of tensors in $\F^2 \otimes \F^3 \otimes \F^3$. There are precisely $20$ $H$-orbits and $17$ $G$-orbits of tensors in $\RR^2 \otimes \RR^3 \otimes \RR^3$.
\end{theorem}

We continue with the line  $\langle x_1,x_2\rangle$ without any assumptions on the field $\F$.
We may assume that $ v_1  =  e_1\otimes e_1 + e_2\otimes e_2 + e_3 \otimes e_3 $, which corresponds to the identity matrix. Recall that the group $H_1$ acts on matrices as $(g,h):X \mapsto gXh$. The characteristic polynomial of a tensor $v$ is $p_v(t) := \det(v_1 - t v)$.


\medskip


Now suppose $\langle v_1,v_2 \rangle$ is equivalent to $\langle v_1,v_3 \rangle$, under an element $(g,h) \in H_1$. Then $g(\alpha v_1 + \beta v_2)h = v_1$, and $gh = \gamma v_1 +\delta v_3$, for some $\alpha,\beta,\gamma,\delta \in \F$. 

If $\beta = 0$, then $\delta = 0$, and $h = \alpha g^{-1}$. Hence $g\langle v_1,v_2\rangle g^{-1} = \langle v_1,v_3\rangle$, and hence $gv_2g^{-1} = \mu v_1 + \nu v_3$ for some $\mu,\nu \in \F$, i.e. $v_2$ is similar to $\mu v_1 + \nu v_2$. Clearly it is also true that if $v_2$ is similar to $\mu v_1 + \nu v_3$, then $\langle v_1,v_2 \rangle$ is equivalent to $\langle v_1,v_3 \rangle$.

If $\beta \ne 0$, then $\delta \ne 0$. Then $h = g^{-1}(\gamma v_1 + \delta v_3)$, and $g(\alpha v_1 + \beta v_2)g^{-1}(\gamma v_1 + \delta v_2) = v_1$. Hence $g(\alpha v_1 + \beta v_2)g^{-1} = (\gamma v_1 + \delta v_3)^{-1}$, and so $(\alpha v_1 + \beta v_2)$ is similar to $(\gamma v_1 + \delta v_3)^{-1}$. It is easy to see that the converse is also true, and hence we have shown the following.
\begin{lemma}
\label{lem:lineequiv}
The line $\langle v_1,v_2 \rangle$ is $H_1$-equivalent to $\langle v_1,v_3 \rangle$ if and only if one of the following hold:
\begin{itemize}
\item
$v_2$ is similar to $\mu v_1 + \nu v_3$ for some $\mu,\nu \in \F$;
\item
$(\alpha v_1 + \beta v_2)$ is similar to $(\gamma v_1 + \delta v_3)^{-1}$ for some $\alpha,\beta,\gamma,\delta \in \F$, $\beta\delta \ne 0$.
\end{itemize}
\end{lemma}

Now two matrices are similar if and only if they have the same rational canonical form. As $v_2$ and $v_3$ both have irreducible characteristic polynomial, their rational canonical forms are the companion matrices of their respective characteristic polynomials.

We recall the definition of the Mobius transformation of a polynomial: if $\phi = \npmatrix{a&b\\c&d}$ is an invertible matrix, and $f$ a polynomial of degree $n$, then 
\[
f^{\phi} (t):= (ct + d)^n f\left(\frac{at+b}{ct+d}  \right)
\]
Note that $f^{a\phi} = f^{\phi}$ for all $a \in \F$, and so we may take $\phi$ to be an element of $\PGL(2,q)$.

We recall also the following straightforward lemma.
\begin{lemma}
Let $f$ be a polynomial of degree $n$, $C(f)$ its companion matrix, and $I$ the identity matrix. Then
\begin{align*}
p_{\alpha I + \beta C(f)}(t) &= \beta^n f \left(\frac{t - \alpha}{\beta}\right) = f^{\phi}(t)\\
p_{C(f)^{-1}}(t) &= -t^n f(0)^{-1} f(-t^{-1}) = f^{\rho}(t),
\end{align*}
where $\phi = \npmatrix{1&-\alpha\\ 0 &\beta}$, $\rho = \npmatrix{0&-1\\ 1 &0}$.
\end{lemma}

This allows us to prove the following lemma.
\begin{lemma}
\label{lem:charpequiv}
Let $v_2$ and $v_3$ be two $n \times n$ matrices with irreducible characteristic polynomials $f$ and $g$ respectively. Then $\langle v_1,v_2 \rangle$ is $H_1$-equivalent to $\langle v_1,v_3 \rangle$ if and only if there exists some $\phi \in \PGL(2,q)$ such that $g = f^{\phi}$.
\end{lemma}

\Proof
We need to consider the two cases from Lemma \ref{lem:lineequiv}. In the first case, $v_2$ is similar to $\mu v_1 + \nu v_3$ if and only if $f(t)= p_{\mu v_1 + \nu v_3}(t) = g^{\phi}(t)$, with $\phi = \npmatrix{1&-\mu\\ 0 &\nu}$.

For the second case, $(\alpha v_1 + \beta v_2)$ is similar to $(\gamma v_1 + \delta v_3)^{-1}$ if and only if 
\[
p_{\alpha v_1 + \beta v_2} (t)= p_{(\gamma v_1 + \delta v_3)^{-1}}(t),
\]
if and only if
\[
 f^{\psi}(t) =  p_{\gamma v_1 + \delta v_3}^{\rho}(t) = g^{\chi\rho}(t),
\]
with $\psi = \npmatrix{1&-\alpha\\ 0 &\beta}$, $\rho = \npmatrix{0&-1\\ 1 &0}$, $\chi = \npmatrix{1&-\gamma\\ 0 &\delta}$, if and only if 
\[
f = g^{\chi\rho\psi^{-1}}.
\]
Now $\chi\rho\psi^{-1} = \frac{1}{\beta}\npmatrix{\gamma \beta&1+\gamma\alpha \\ -\delta \beta &\delta \alpha}$, and it is clear that every element of $\PGL(2,q)$ can be represented by one of the forms $\npmatrix{1&-\mu\\ 0 &\nu}$, $\frac{1}{\beta}\npmatrix{\gamma \beta&1+\gamma\alpha \\ -\delta \beta &\delta \alpha}$, $\nu,\beta,\delta \ne 0$, proving the claim. \qed

\begin{theorem}
\label{thm:rk3}
There is precisely one orbit of tensors in $\Fq^2 \otimes \Fq^3 \otimes \Fq^3$ whose first contraction space is a line without singular points in $\PG(\Fq^3 \otimes \Fq^3)$.
\end{theorem}

\Proof
It suffices by Lemma \ref{lem:charpequiv} to show that every two monic irreducible polynomials of degree $3$ are related by a Mobius transformation. Let $f$ be a monic irreducible of degree $3$. Then the stabiliser of $f$ under $\PGL(2,q)$ has order $3$, and hence the size of the orbit of $f$ under this action is $\frac{|\PGL(2,q)|}{3} = \frac{q^3-q}{3}$, which is precisely the number of monic irreducible polynomials of degree $3$, proving the claim. \qed

Note that there may be an infinite number of orbits for some fields, such as $\QQ$. 

\begin{remark}
We remark that in $\PG(\Fq^n \otimes \Fq^n)$ with $n > 3$, there is more than one orbit of lines of constant rank $n$ over the finite field $\Fq$. For example if $n$ is prime, a similar argument shows that there are at least $\frac{q^{n-1}-1}{q^2-1}$ orbits, arising from lines $\langle v_1,v_2\rangle$ where $v_2$ has irreducible characteristic polynomial. 
\end{remark}

\subsection{Canonical forms of tensors in $\F^2\otimes \F^3 \otimes \F^3$}\label{sec:summary}
In this section we list the canonical forms for the representatives of the $G$-orbits in $\F^2\otimes \F^3 \otimes \F^3$. The special cases of Theorems \ref{thm:233} and \ref{thm:223} for the field $\F=\F_2$ were computed in \cite{BrHu2012}.

The previous subsections are summarized in the following theorem.

\begin{theorem}
\label{thm:233}
If $\F$ is a finite field, then there are precisely $21$ $H$-orbits and $18$ $G$-orbits of tensors in $\F^2 \otimes \F^3 \otimes \F^3$. For any algebraically closed field $\F$, there are precisely $18$ $H$-orbits and $15$ $G$-orbits of tensors in $\F^2 \otimes \F^3 \otimes \F^3$. There are precisely $20$ $H$-orbits and $17$ $G$-orbits of tensors in $\RR^2 \otimes \RR^3 \otimes \RR^3$.
\end{theorem}

\Proof Follows from Theorems \ref{thm:main} and \ref{thm:rk3}. \qed

The canonical forms for the 18 $G$-orbits in the finite field case can be extracted from the proof of the classification of the orbits. If $e_1,e_2,e_3$ is a basis for $\F^3$ and $e=e_1 \otimes e_1+ e_2\otimes e_2 + e_3 \otimes e_3$, then the canonical forms are as follows. The third column contains the rank distribution $r_1(A)$ of the first contraction space of a representative $A$.
\medskip

\begin{tabular}{|l|l|l|}
\hline
Orbit&Canonical form&$r_1(A)$\\
\hline
$o_0$& 0&$[0,0,0]$ \\
 $o_1$& $ e_1 \otimes e_1 \otimes e_1 $&$[1,0,0]$\\
$o_2$& $ e_1 \otimes (e_1 \otimes e_1+ e_2\otimes e_2) $&$[0,1,0]$\\
$o_3$& $ e_1 \otimes e $&$[0,0,1]$\\
$o_4$& $ e_1 \otimes e_1 \otimes e_1+ e_2\otimes e_1 \otimes e_2  $&$[q+1,0,0]$\\
 $o_5$& $ e_1 \otimes e_1 \otimes e_1+ e_2\otimes e_2 \otimes e_2  $&$[2,q-1,0]$\\
 $o_6$ & $ e_1 \otimes e_1 \otimes e_1+ e_2\otimes (e_1 \otimes e_2 + e_2 \otimes e_1)  $&$[1,q,0]$\\
 $o_7$& $ e_1 \otimes e_1 \otimes e_3+ e_2\otimes (e_1 \otimes e_1 + e_2 \otimes e_2)  $&$[1,q,0]$\\
 $o_8$& $ e_1 \otimes e_1 \otimes e_1+ e_2\otimes (e_2 \otimes e_2 + e_3 \otimes e_3)  $&$[1,1,q-1]$\\
$o_9$& $ e_1 \otimes e_3 \otimes e_1+ e_2\otimes e  $&$[1,0,q]$\\
$o_{10}$&  $ e_1\otimes (e_1\otimes e_1+ e_2\otimes e_2+u e_1\otimes e_2) +  
e_2\otimes (e_1\otimes e_2+v e_2\otimes e_1),$&$[0,q+1,0]$\\
 & $v\lambda^2+uv\lambda - 1 \neq 0$ for all
$\lambda \in \F$& \\
$o_{11}$& $ e_1\otimes (e_1 \otimes e_1 + e_2 \otimes e_2)+ 
e_2\otimes (e_1 \otimes e_2 + e_2 \otimes e_3)$&$[0,q+1,0]$\\
$o_{12}$& $ e_1\otimes (e_1 \otimes e_1 + e_2 \otimes e_2)+ 
e_2\otimes (e_1 \otimes e_3 + e_3 \otimes e_2)$&$[0,q+1,0]$\\
$o_{13}$& $ e_1\otimes (e_1 \otimes e_1 + e_2 \otimes e_2)+ 
e_2\otimes (e_1 \otimes e_2 + e_3 \otimes e_3)$&$[0,2,q-1]$\\
$o_{14}$& $ e_1\otimes (e_1 \otimes e_1 + e_2 \otimes e_2)+ e_2\otimes (e_2 \otimes e_2 + e_3 \otimes e_3)$&[0,3,q-2]\\
$o_{15}$&  $ e_1\otimes (e+u e_1\otimes e_2) + 
e_2\otimes (e_1\otimes e_2+v e_2\otimes e_1),$&$[0,1,q]$\\
 & $v\lambda^2+uv\lambda - 1 \neq 0$ for all
$\lambda \in \F$& \\
$o_{16}$& $ e_1\otimes e+ e_2\otimes (e_1 \otimes e_2 + e_2 \otimes e_3)$&$[0,1,q]$\\
$o_{17}$& $ e_1\otimes e+ 
e_2\otimes (e_1\otimes e_2 + e_2\otimes  e_3 + e_3\otimes (\alpha e_1 + \beta e_2 + \gamma e_3)),$ &$[0,0,q+1]$\\
 & $\lambda^3+\gamma \lambda^2- \beta \lambda+ \alpha \neq 0$ for all $\lambda \in \F$ & \\
 \hline
\end{tabular}
\medskip

This classification also allows us to classify tensors in $\F^2 \otimes \F^2 \otimes \F^3$. Orbits $o_0,o_1,o_2,o_4,o_4^T,o_5,o_6$, and $o_{10}$ are contained in a space isomorphic to $\F^2 \otimes \F^2 \otimes \F^2$. The orbits corresponding to $o_2$, $o_4$ and $o_4^T$ are $G$-equivalent in $\F^2 \otimes \F^2 \otimes \F^2$. 

Furthermore, orbits $o_7$ and $o_{11}$ are contained in a space isomorphic to $\F^2 \otimes \F^2 \otimes \F^3$. The orbits corresponding to $o_2$ and $o_4$ are $G$-equivalent in $\F^2 \otimes \F^2 \otimes \F^3$, but not $G$-equivalent to $o_4^T$. This implies the following.

\begin{theorem}
\label{thm:223}
If $\F$ is a finite or real field, then there are precisely $10$ $H$-orbits and $9$ $G$-orbits of tensors in $\F^2 \otimes \F^2 \otimes \F^3$. For any algebraically closed field $\F$, there are precisely $9$ $H$-orbits and $8$ $G$-orbits of tensors in $\F^2 \otimes \F^2 \otimes \F^3$. 
\end{theorem}

\begin{remark}
In \cite[Table 4]{Nurmiev2000}, representatives of each $G$-orbit in $\CC^3 \otimes \CC^3 \otimes \CC^3$ are given. 
For the reader's convenience, we include the correspondence between the orbits $o_0, \ldots o_{17}$ and the orbits as listed in \cite{Nurmiev2000}. Note that $o_2$ and $o_4$ are equivalent under $G$ here, while $o_{10},o_{15}$ and $o_{17}$ are empty.
%
\[
\begin{array}{|c|c|c|c|c|c|c|c|c|c|c|c|c|c|c|c|}
\hline
&o_0 &o_1 & o_2 &o_3&o_4&o_5&o_6&o_7&o_8& o_9&o_{11}&o_{12}&o_{13}&o_{14}&o_{16}\\
\hline
\cite{Nurmiev2000}&25&24&23&22&23&20&21&19&15&16&18&17&12&9&13\\
\hline
\cite{Thrall1938}&&&&&&&&&&&&&&&\\
\hline
\end{array}
\]
\end{remark}

\end{document}